\documentclass[12pt]{article}
\usepackage{amsmath,amssymb,amsthm}

\usepackage{fullpage}

\title{Ramanujan's Perimeter of an Ellipse}

\author{MARK B. VILLARINO \\[12pt]
        Escuela de Matem\'atica,\\
        Universidad de Costa Rica,\\
        2060 San Jos\'e, Costa Rica}

\date{\today}

\newtheorem{theorem}{THEOREM}

\newtheorem{Cor}{COROLLARY}
\theoremstyle{definition}

\newtheorem{claim}{CLAIM}
\newtheorem{corollary}{COROLLARY}

\numberwithin{equation}{section}

\newcommand{\dis}{\displaystyle} 

\begin{document}

\maketitle

\begin{abstract}
We present a detailed analysis of \textsc{Ramanujan}'s most accurate approximation to the perimeter of an ellipse. 
\end{abstract}

\tableofcontents
\section{Introduction}

Let $a$ and $b$ be the semi-major and semi-minor axes of an ellipse with perimeter $p$ and whose eccentricity is $k$.  The final sentence of \textsc{Ramanujan}'s famous paper \emph{Modular Equations and Approximations to $\pi$,} \cite{R}, says:\emph{\begin{quotation}``  The following approximation for  $p$ [was] obtained empirically:\begin{equation}\fbox{$\dis
p=\pi\left\{(a+b)+\frac{3(a-b)^{2}}{10(a+b)+\sqrt{a^{2}+14ab+b^{2}}}+\varepsilon\right\}$}
\end{equation}\indent where  $\varepsilon$ is about $\dfrac{3ak^{20}}{68719476736}.$"\end{quotation}}\textsc{Ramanujan} never explained his ``empirical" method of obtaining this approximation, nor ever subsequently returned to this approximation, neither in his published work, nor in his Notebooks \cite{
B}.  Indeed, although the Notebooks does contain the above approximation (see Entry 3 of Chapter XVIII) the statement there does not even mention his asymptotic error estimate stated above.

Twenty years later \textsc{Watson} \cite{W} claimed to have proven that \textsc{Ramanujan}'s approximation is \textbf{\emph{in defect}}, but he never published his proof.

In 1978, we established the following optimal version of \textsc{Ramanujan}'s approximation:
\begin{theorem}(\textbf{Ramanujan's Approximation Theorem}) \textsc{Ramanujan}'s approximative perimeter
\begin{equation}\fbox{$\dis
p_{R}:=\pi\left\{(a+b)+\frac{3(a-b)^{2}}{10(a+b)+\sqrt{a^{2}+14ab+b^{2}}}\right\}$}
\end{equation}
\textbf{underestimates} the \textbf{true} perimeter, $p$, by\begin{equation}
\fbox{$\dis\epsilon:=\pi(a+b)\cdot\theta(\lambda)\cdot\lambda^{10},
$}\end{equation}where\begin{equation}
\lambda:=\frac{a-b}{a+b},
\end{equation}and where the function $\theta(\lambda)$ \textbf{grows monotonically} in $0\leqslant \lambda\leqslant 1$ while at the same time it satisfies the \textbf{optimal} inequalities
\begin{equation}\fbox{$\dis
\frac{3}{2^{17}}<\theta(\lambda)\leqslant \frac{14}{11}\left(\dfrac{22}{7}-\pi\right)$}
\end{equation}
\end{theorem}\hfill$\Box$

Please take note of the striking form of the sharp upper bound since it involves the number $\left(\dfrac{22}{7}-\pi\right)$ which measures \emph{the accuracy of \textsc{Archimedes}' famous approximation, $\dfrac{22}{7},$ to the transcendental number $\pi$!}

\begin{Cor}The error in defect, $\epsilon$, as a function of $\lambda$, \textbf{grows monotonically} for $0\leqslant \lambda\leqslant 1.$
\end{Cor}\hfill$\Box$

\begin{Cor}The error in defect, $\epsilon$, as a function of the \textbf{eccentricity}, $e$, is given by \begin{equation}
\fbox{$\dis \epsilon(e):=a\left\{\delta(e)\left(\frac{2}{1+\sqrt{1-e^{2}}}\right)^{19}\right\} e^{20}$}.
\end{equation}Moreover, $\epsilon(e)$ \textbf{grows monotonically} with $e$, $0\leqslant e\leqslant 1$, while $\delta(e)$ satisfies the \textbf{optimal inequalities}\begin{equation}
\fbox{$\dis\frac{3\pi}{68719476736}<\delta(e)\leqslant \dfrac{\dfrac{7}{11}\left(\dfrac{22}{7}-\pi\right)}{2^{18}}$}
\end{equation} 
\end{Cor}\hfill$\Box$
\\

This \textbf{Corollary 2} explains the significance of \textsc{Ramanujan}'s own error estimate in (1.1).  The latter is an asymptotic \emph{lower bound} for $\epsilon(e)$ but it is not the optimal one.  That is given in (1.7).

\section{Later History}

We sent an (updated) copy of our 1978 preprint to Professor \textsc{Bruce Berndt} in 1988 and he subsequently quoted its conclusions in his edition of Volume 3 of the Notebooks (see p. 150 \cite{B}). However the details of our proofs have never been published and so we have decided to present them in this paper.

\textsc{Berndt}'s discussion of \textsc{Ramanujan}'s approximation includes \textsc{Almkvist}'s very plausible suggestion that \textsc{Ramanujan}'s ``empirical process" was to develop a\emph{ continued fraction expansion} of \textsc{Ivory}'s infinite series for the perimeter (\cite{A}) as well as a proof, due independently to \textsc{Almkvist} and \textsc{Askey}, of our fundamental lemma (see \S 3).  However, their proof is different from ours.

The most recent work on the subject has been carried out by \textsc{R. Barnard, K. Pearce,} and \textsc{K. Richards} in \cite{BPR} and was published in the year $2000$.  They also prove the major conclusion in our fundamental lemma, but their methods too are quite different from ours.

\section{Fundamental Lemma}

\begin{theorem}(\textbf{Fundamental Lemma}) Define the functions $\mathbf{A}(x)$ and $\mathbf{B}(x)$ and the coefficients $A_{n}$ and $B_{n}$ by:\begin{align}
\mathbf{A}(x)&:=1+\frac{3x}{10+\sqrt{4-3x}}:=1+A_{1}x+A_{2}x^{2}+\cdots   \\
 \mathbf{B}(x)&:=\sum_{n=0}^{\infty}\left\{\frac{1}{2n-1}\frac{1}{4^{n}}\binom{2n}{n}\right\}^{2}x^{n}:=1+B_{1}x+B_{2}x^{2}+\cdots.\end{align}
 Then:\begin{align}
    &A_{1}=B_{1},  \ A_{2}=B_{2}, \ A_{3}=B_{3}, \ A_{4}=B_{4}   \\
    &A_{5}<B_{5}, \ A_{6}<B_{6}, \cdots, A_{n}<B_{n}, \cdots
    \end{align}      
where the strict inequalities in (3.4) are valid for all $n\geqslant 5.$
\end{theorem}
\begin{proof}

\textbf{First we prove} $\mathbf{(3.3)}$.  We read this off directly from the numerical values of the expansion:\begin{align*}
 A_{1} &= B_{1}=\frac{1}{4}  \\
    A_{2} &= B_{2}=\frac{1}{16}  \\
    A_{3} &= B_{3}=\frac{1}{64}  \\
    A_{4} &= B_{4}=\frac{25}{4096}.    
\end{align*}
\\
\textbf{Now we prove} $\mathbf{(3.4)}$  For $A_{5}$, $B_{5}$, $A_{6}$, and $B_{6}$ we verify $(3.4)$ directly from their explicit numerical values.  Namely,\begin{align*}
 A_{5} &=\frac{47\frac{1}{2}}{2^{14}}, \ B_{5}=\frac{49}{2^{14}},\ \Rightarrow A_{5}-B_{5}=\frac{-\frac{3}{2}}{2^{14}} <0  \\
     A_{6} &=\frac{803}{2^{21}}, \ B_{6}=\frac{882}{2^{21}},\ \Rightarrow A_{6}-B_{6}=\frac{-79}{2^{21}} <0. 
\end{align*}\textbf{\emph{Therefore it is sufficient to prove }}\begin{equation}
A_{n}<B_{n}
\end{equation}
\textbf{\emph{for all}}\begin{equation}
n\geqslant 7.
\end{equation}

Now the \textbf{\emph{explicit}} formula for $A_{n}$ is\begin{equation}
A_{n}=a_{n-1}+a_{n-2}+a_{n-3}+\cdots+a_{1}+a_{0}
\end{equation}where\begin{equation}
\fbox{$\begin{array}{rll}
     a_{n-1} &:=\dis\frac{1}{2n-3}\frac{1}{16^{n}}\binom{2n-2}{n-1}3^{n-1}    \\
      a_{n-2} &:=\dis\frac{1}{2n-5}\frac{1}{16^{n-1}}\binom{2n-2}{n-1}3^{n-2} \left(\frac{-1}{2^{5}}\right)\\
    .\ \ \ &\ \ \ \ \ \  .\\
      .\ \ \ &\ \ \ \ \ \  .\\
      .\ \ \ &\ \ \ \ \ \  .\\
     a_{1} &:=\dis\frac{1}{2\cdot 1-1}\frac{1}{16^{2}}\binom{2}{1}3^{n-2} \left(\frac{-1}{2^{5}}\right)^{n-2}\\ 
     a_{0}&:=\dis\frac{4}{16} \left(\frac{-1}{2^{5}}\right)^{n-1}.  
\end{array}$}
\end{equation}Next we write\begin{equation}
A_{n}=a_{n-1}\left(1+\frac{a_{n-2}}{a_{n-1}}+\frac{a_{n-3}}{a_{n-1}}+\frac{a_{n-4}}{a_{n-1}}+\cdots+\frac{a_{1}}{a_{n-1}}+\frac{a_{0}}{a_{n-1}}\right)
\end{equation}and assert:
\begin{claim}\indent\\

\emph{The ratios $\dfrac{a_{n-k-1}}{a_{n-k}}$ \textbf{decrease monotonically} in absolute value as $k$ increases from $k=1$ to $k=n-1.$
}

\end{claim}
\begin{proof}
For $k=1, \cdots, n-2,$\begin{align*}
\left|\dfrac{a_{n-k-1}}{a_{n-k}}\right|    &=\left(1+\frac{2}{2n-2k-3}\right)\left(1+\frac{1}{4n-4k-2}\right)\frac{1}{12}   \\
    & \leqslant \frac{1}{6}\ \ (\text{which is the worst case and occurs when $k=n-2$)}\\
    &<1
\end{align*}For $k=n-1$, $$\left|\frac{a_{0}}{a_{1}}\right|=\frac{1}{3}<1.$$This completes the proof.
\end{proof}
\begin{claim}\indent\\

\emph{The ratios $\dfrac{a_{n-k-1}}{a_{n-k}}$\textbf{alternate in sign.}}
\end{claim}
\begin{proof}
This is a consequence of the definition of the $a_{k}.$

\end{proof}
By \textbf{CLAIM 1.} and \textbf{CLAIM 2.} we can write $(3.9)$ in the form \begin{align*}
   A_{n} &=a_{n-1}(1-\text{something positive and smaller than $1$})   \\
    &<a_{n-1}.  
\end{align*}Therefore, to prove $(3.8)$ for $n\geqslant 7$, \textbf{\emph{it suffices to prove}}\begin{equation}
a_{n-1}<B_{n}
\end{equation}\textbf{\emph{for all}} $n\geqslant 7$.

By $(3.8)$ and the definition of $B_{n}$, this last afirmation is equivalent to proving$$\frac{1}{2n-3}\frac{1}{16^{n}}\binom{2n-2}{n-1}3^{n-1}<\left\{\frac{1}{2n-1}\frac{1}{4^{n}}\binom{2n}{n}\right\}^{2},$$ which, after some algebra, \textbf{\emph{reduces to proving the implication}} $$n\geqslant 7 \Rightarrow\frac{\dis\frac{n}{2}\cdot\frac{2n-1}{2n-3}}{\dis\binom{2n}{n}}\cdot3^{n-1}<1.$$
If we define for all integers $n\geqslant 7$\begin{equation}
f(n):=\frac{\dis\frac{n}{2}\cdot\frac{2n-1}{2n-3}}{\dis\binom{2n}{n}}\cdot3^{n-1}
\end{equation}then the affirmation $(3.10)$ turns out to be\textbf{\emph{ equivalent to}}\begin{equation}
n\geqslant 7\Rightarrow f(n)<1
\end{equation}

This latter affirmation is a consequence of the following two conditions:\begin{enumerate}
  \item $f(7)<1.$
  \item $f(7)>f(8)>f(9)>\cdots>f(k)>f(k+1)>\cdots$
  \end{enumerate}

\textbf{Proof of 1.}  By direct numerical computation $$f(7)=\frac{1701}{1936}<1$$
\hfill$\Box$

\textbf{Proof of 2.}  \textbf{We must show}$$k\geqslant 7 \Rightarrow f(k)>f(k+1).$$If we define\begin{equation}g(k):=\frac{f(k)}{f(k+1)},\end{equation}then \textbf{we must show}\begin{equation}
k\geqslant 7\Rightarrow g(k)>1.
\end{equation}Using the definition $(3.11)$ of $f(n)$ and the definition $(3.14)$ of $g(n)$, and reducing algebraically we find $$g(k)=\frac{2k}{6k-9}\left(\frac{2k-1}{k+1}\right)^{2},$$ and \textbf{we must show} that\begin{equation}
k\geqslant 7 \Rightarrow \frac{2k}{6k-9}\left(\frac{2k-1}{k+1}\right)^{2}>1.
\end{equation}Define the rational function of the real variable $x$:\begin{equation}
g(x):=\frac{2x}{6x-9}\left(\frac{2x-1}{x+1}\right)^{2}.
\end{equation}Then the graph of $y=g(x)$ has a vertical asymptote at $x=\frac{3}{2}$ and\begin{equation}
\lim_{x\rightarrow\frac{3}{2}^{+}}g(x)=+\infty.
\end{equation}Moreover, the derivative of $g(x)$ is given by:$$g'(x)=\frac{2(2x^{2}-7x+1)}{x(x+1)(2x-1)(2x+3)}$$which implies that $$g'(x)\begin{cases}
     <0 & \text{if $\frac{3}{2}<x<\frac{7+\sqrt{41}}{4}$ }, \\
       =0&\text{if $x=\frac{7+\sqrt{41}}{4}$},\\
      >0& \text{if $x>\frac{7+\sqrt{41}}{4}$}.
\end{cases}$$Therefore, for $x\geqslant \frac{3}{2}$ g(x) \textbf{\emph{decreases}} from ``$+\infty$" at $x=\frac{3}{2}$ (see $(3.17)$) to an \textbf{\emph{absolute minimum value}} (in $\frac{3}{2}\leqslant x<\infty$)
$$g\left(\frac{7+\sqrt{41}}{4}\right)=1+\frac{37-\sqrt{41}}{399+69\sqrt{41}}=1.0363895208\cdots$$and then \textbf{\emph{increases monotonically}} as $x\rightarrow \infty$ to its asymptotic limit $y=\frac{4}{3}.$Therefore\begin{align*}
g(x)    &>1.03638\cdots \text{for all $x>\frac{3}{2}, x\neq\frac{7+\sqrt{41}}{4}$}   \\
 \Rightarrow g(x)   &>1.03638\cdots\  \text{for all integers $n\geqslant 2.$}  \\
\Rightarrow f(n)&>(1.03638\cdots)f(n+1)\  \text{for all integers $n\geqslant 2.$}\\
\Rightarrow f(n)&>f(n+1)\  \text{for all integers $n\geqslant 2.$}\\
\Rightarrow f(7)&>f(8)>f(9)>\cdots\\
\end{align*}which implies that the condition 2. holds.  Moreover we conclude that\begin{align*}
  & f(n) <1 \ \text{for all integers $n\geqslant 7$}   \\
   \Rightarrow &(3.10) \text{holds for all integers $n\geqslant 7$}\\
   \Rightarrow &(3.5)\  \text{holds for all integers $n\geqslant 7$}\\
   \Rightarrow &(3.4)\  \text{holds for all integers $n\geqslant 5$}\\  
\end{align*}and this completes the proof of the Fundamental Lemma.
\end{proof}
\section{Ivory's Identity}

In $1796$, \textsc{J. Ivory} \cite{I} published the following identity (in somewhat different notation):

\begin{theorem}(\textbf{Ivory's Identity}) If $0\leqslant x\leqslant 1$ then the following formula for $\mathbf{B}(x)$ is valid:\begin{equation}
\fbox{$\dis \frac{1}{\pi}\int_{0}^{\pi}\sqrt{1+2\sqrt{x}\cos(2\phi)+x}~d\phi=\sum_{n=0}^{\infty}\left\{\frac{1}{2n-1}\frac{1}{4^{n}}\binom{2n}{n}\right\}^{2}x^{n}\equiv\mathbf{B}(x)$}
\end{equation}
\end{theorem}

\begin{proof}

We sketch his elegant proof.  \begin{align*}
  &\frac{1}{\pi}\int_{0}^{\pi}\sqrt{1+2\sqrt{x}\cos(2\phi)+x}~d\phi =\frac{1}{\pi}\int_{0}^{\pi}\sqrt{1+\sqrt{x}e^{2i\phi}}\sqrt{1+\sqrt{x}e^{-2i\phi}}~d\phi&  \\
    & = \frac{1}{\pi}\int_{0}^{\pi}\sum_{m=0}^{\infty}\left\{\frac{1}{2m-1}\frac{1}{4^{m}}\binom{2m}{m}(\sqrt{x})^{m}e^{2\pi i m \phi}\right\}\sum_{n=0}^{\infty}\left\{\frac{1}{2n-1}\frac{1}{4^{n}}\binom{2n}{n}(\sqrt{x})^{n}e^{-2\pi i n \phi}\right\}~d\phi
\\
&=\frac{1}{\pi}\sum_{m=0}^{\infty}\left\{\frac{1}{2m-1}\frac{1}{4^{m}}\binom{2m}{m}(\sqrt{x})^{m}\right\}\sum_{n=0}^{\infty}\left\{\frac{1}{2n-1}\frac{1}{4^{n}}\binom{2n}{n}(\sqrt{x})^{n}\right\}\int_{0}^{\pi}e^{2\pi i(m-n)\phi}~d\phi\\
&=\sum_{n=0}^{\infty}\left\{\frac{1}{2n-1}\frac{1}{4^{n}}\binom{2n}{n}\right\}^{2}x^{n}
\end{align*}
\end{proof}

We will need the following evaluation in our investigation of the accuracy of \textsc{Ramanujan}'s approximation.
\begin{corollary}\begin{equation}
\fbox{$\dis\mathbf{B}(1)=\frac{4}{\pi}$}
\end{equation}
\end{corollary}
\begin{proof}

By \textsc{Ivory}'s identity,\begin{align*}
  \mathbf{B}(1)=  &\frac{1}{\pi}\int_{0}^{\pi}\sqrt{1+2\sqrt{1}\cos(2\phi)+1}~d\phi   \\
    &=\frac{1}{\pi}\int_{0}^{\pi}\sqrt{2+2\cos(2\phi)}~d\phi\\
    &=\frac{1}{\pi}\int_{0}^{\pi}\sqrt{4\cos^{2}(\phi)}~d\phi\\
    &= \frac{4}{\pi} 
\end{align*}

\end{proof}

\section{The Accuracy Lemma}

\begin{theorem}(\textbf{Accuracy Lemma}) For $0\leqslant x\leqslant 1$, the function
\begin{equation}
\mathbf{A}(x):=1+\frac{3x}{10+\sqrt{4-3x}}
\end{equation}
\textbf{underestimates} the function
\begin{equation}
\mathbf{B}(x):=\sum_{n=0}^{\infty}\left\{\frac{1}{2n-1}\frac{1}{4^{n}}\binom{2n}{n}\right\}^{2}x^{n}
\end{equation}by a \textbf{discrepancy}, $\Delta(x)$ which is \textbf{never more} than $\dis\left(\frac{4}{\pi}-\frac{14}{11}\right)x^{5}$ and which is \textbf{always more} than $\dfrac{3}{2^{17}}x^{5}$:\begin{equation}
\fbox{$\dis\dfrac{3}{2^{17}}x^{5}<\Delta(x)\leqslant \dis\left(\frac{4}{\pi}-\frac{14}{11}\right)x^{5}$}
\end{equation} Moreover, the constants $\dis\left(\frac{4}{\pi}-\frac{14}{11}\right)$ and $\dfrac{3}{2^{17}}x^{5}$ are \textbf{the best possible.}

\end{theorem}

\begin{proof}
By the definition of $\mathbf{A}(x)$ and $\mathbf{B}(x)$ given in \textbf{Theorem 1.}, the discrepancy $\Delta(x)$ is given by the series\begin{align*}
\Delta(x)    &:=\mathbf{B}(x)-\mathbf{A}(x)   \\
    & =(B_{5}-A_{5})x^{5}+(B_{6}-A_{6})x^{6}+\cdots \\
    &:=\delta_{5}x^{5}+\delta_{6}x^{6}+\cdots,
\end{align*}
\noindent where, again by \textbf{Theorem 1.}, $$\delta_{k}>0 \ \ \ \text{for $k=5, 6, \cdots.$}$$

On the one hand \begin{align*}
   \Delta(x) &=x^{5}(\delta_{5}+\delta_{6}x+\cdots)   \\
    & \leqslant x^{5}(\delta_{5}+\delta_{6}+\delta_{7}+\cdots)\\
    &=x^{5}\Delta(1)\\
    &=x^{5}\{\mathbf{B}(1)-\mathbf{A}(1)\}\\
    &=x^{5}\left(\frac{4}{\pi}-\frac{14}{11}\right)
\end{align*}where we used \textbf{Corollary 1} of \textsc{Ivory}'s identity in the last equality.  Therefore\begin{center}
\boxed{\Delta(x)\leqslant \dis\left(\frac{4}{\pi}-\frac{14}{11}\right)x^{5}.}
\end{center}This is half of the accuracy lemma.  Moreover the constant $\dis\left(\frac{4}{\pi}-\frac{14}{11}\right)$ is \textbf{\emph{assumed}} for $x=1$ and thus cannot be replaced by anything smaller, i.e., it is \textbf{\emph{the best possible}} constant.

On the other hand, we can write$$\Delta(x)=x^{5}\{\delta_{5}+G(x)\},$$where$$G(x):=\delta_{6}x+\delta_{7}x^{2}+\cdots\Rightarrow\begin{cases}
   G(x)\geqslant 0   & \text{for all $0\leqslant x\leqslant 1$ }, \\
    G(x)\rightarrow 0  & \text{as $x\rightarrow 0$}.
\end{cases}$$This shows that \begin{center}
\boxed{\Delta(x)>\delta_{5}x^{5}=\dfrac{3}{2^{17}}x^{5}}
\end{center}and that $$\lim_{x\rightarrow 0}\frac{\Delta(x)}{x^{5}}=\frac{3}{2^{17}}.$$ This proves both the other inequality in the theorem and the 
  \textbf{\emph{optimality}} of the constant $\delta_{5}=\dfrac{3}{2^{17}},$ i.e., that it cannot be replaced by any larger constant.
\\

This completes the proof of the \textbf{Accuracy Lemma}.\end{proof}
\section{The Accuracy of Ramanujan's Approximation}
Now we can achieve the main goal of this paper, namely to prove \textbf{\emph{Ramanujan's Approximation Theorem}}.

First we express the perimeter of an ellipse and \textsc{Ramanujan}'s approximative perimeter in terms of the functions $\mathbf{A}(x)$ and $\mathbf{B}(x)$.
\begin{theorem} If $p$ is the perimeter of an ellipse with semimajor axes $a$ and $b$, and if $p_{R}$ is \textsc{Ramanujan}'s approximative perimeter, then:
\begin{equation}
\fbox{$\dis \begin{array}{rll }
     p&=\pi(a+b)\cdot\mathbf{B}\left\{\left(\dfrac{a-b}{a+b}\right)^{2}\right\}\\
     &\\
      p_{R}&=\pi(a+b)\cdot\mathbf{A}\left\{\left(\dfrac{a-b}{a+b}\right)^{2}\right\}.   
\end{array}$}      
\end{equation}
\end{theorem}
\begin{proof}
We begin with \textbf{\emph{Ivory's Identity}} (\S4) and in it we substitute $x:=\left(\dfrac{a-b}{a+b}\right)^{2}.$  Then the integral becomes\begin{align*}
\frac{1}{\pi}\int_{0}^{\pi}\sqrt{1+2\sqrt{\left(\dfrac{a-b}{a+b}\right)^{2}}\cos(2\phi)+\left(\dfrac{a-b}{a+b}\right)^{2}}~d\phi    &= \frac{4}{\pi(a+b)}\int_{0}^{\frac{\pi}{2}}(a^{2}\sin^{2}\phi+b^{2}\cos^{2}\phi)  ~d\phi
     \end{align*}and therefore\begin{align*}
\mathbf{B}\left\{\left(\dfrac{a-b}{a+b}\right)^{2}\right\}    &= \frac{4}{\pi(a+b)}\int_{0}^{\frac{\pi}{2}}(a^{2}\sin^{2}\phi+b^{2}\cos^{2}\phi)  ~d\phi   \\
      \end{align*}But, it is well known (\textsc{Berndt} \cite{B}) that \textbf{\emph{the perimeter, $p$, of an ellipse with semiaxes $a$ and $b$ is given by }}$$p=4\int_{0}^{\frac{\pi}{2}}(a^{2}\sin^{2}\phi+b^{2}\cos^{2}\phi)  ~d\phi,$$ and thus\begin{equation}
p=\pi(a+b)\cdot\mathbf{B}\left\{\left(\dfrac{a-b}{a+b}\right)^{2}\right\}.
\end{equation}

Moreover, some algebra shows us that \begin{align*}
\mathbf{A}\left\{\left(\dfrac{a-b}{a+b}\right)^{2}\right\}&=1+\frac{3\left(\dfrac{a-b}{a+b}\right)^{2}}{10+\sqrt{4-3\left(\dfrac{a-b}{a+b}\right)^{2}}}\\
&=\frac{1}{a+b}\left\{(a+b)+\frac{3(a-b)^{2}}{10(a+b)+\sqrt{a^{2}+14ab+b^{2}}}\right\}
\end{align*}and we conclude that \textbf{\textsc{Ramanujan}'s \emph{approximative formula, $p_{R}$ is given by}}\begin{equation}
p_{R}=\pi(a+b)\mathbf{A}\left\{\left(\dfrac{a-b}{a+b}\right)^{2}\right\}.
\end{equation}
\end{proof}

The formula for $p$ above was the object of \textsc{Ivory}'s original paper \cite{I}.

Now we complete the proof of \textbf{Theorem 1.}

\begin{proof}
Writing $$\lambda:=\frac{a-b}{a+b},$$ and using the notation of the statement of \textbf{Theorem 1.} we conclude that\begin{align*}
\epsilon&:=\pi(a+b)\cdot\theta(\lambda)\cdot\lambda^{10}      \\
    &= \pi(a+b)\cdot\frac{\Delta(\lambda^{2})}{\lambda^{10}}\cdot\lambda^{10} 
\end{align*}where\begin{equation}
\theta(\lambda)\equiv\frac{\Delta(\lambda^{2})}{\lambda^{10}}=\delta_{5}+\delta_{6}\lambda^{2}+\cdots.
\end{equation}

Now we apply the \textbf{\emph{Accuracy Lemma}} and the proof is complete.

\end{proof}


\end{document}